\documentstyle[11pt]{article}
\newfont{\bbb}{msbm10 scaled \magstep1}

\newcommand{\Z}{\mbox{\bbb Z}}

\newcommand{\beqa}{\begin{eqnarray*}}
\newcommand{\eeqa}{\end{eqnarray*}}
\newcommand{\beqn}{\begin{eqnarray}}
\newcommand{\eeqn}{\end{eqnarray}}
\newcommand{\eeq}{\end{equation}}
\newcommand{\beq}{\begin{equation}}
\newcommand{\bdes}{\begin{description}}
\newcommand{\edes}{\end{description}}
\newcommand{\ra}{\rightarrow}

\newcommand{\be}{\begin{equation}}
\newcommand{\ee}{\end{equation}}

\newtheorem{defn}{Definition}

\newtheorem{thm}{Theorem}
\newtheorem{prop}{Proposition}
\newtheorem{lem}{Lemma}
\newtheorem{rem}{Remark}

\newcommand{\bed}{\begin{defn}}
\newcommand{\bet}{\begin{thm}}
\newcommand{\bep}{\begin{prop}}
\newcommand{\bel}{\begin{lem}}
\newcommand{\brk}{\begin{rem}}

\newcommand{\eed}{\end{defn}}
\newcommand{\eet}{\end{thm}}
\newcommand{\eep}{\end{prop}}
\newcommand{\erk}{\end{rem}}
\newcommand{\eecor}{\end{cor}}

\newcommand{\non}{\noindent}
\newcommand{\no}{\noindent}

 \newcommand{\proofhead}[1]{\par\pagebreak[1]\noindent{\bf#1.\ }}
 \newcommand{\pf}{\proofhead{Proof}}

\newcommand{\qed}{{\unskip\nolinebreak[1]\hspace{1.5em}\mbox{}\nolinebreak
    \hfill$\Box$\parfillskip=0pt\finalhyphendemerits=0\par\pagebreak[1]}}

\begin{document}
\title{ \bf A generalization of Scheunert's Theorem on cocycle twisting
 of color Lie algebras} 
\author{ Horia C. Pop \\ Department of Mathematics, University of Iowa\\
Iowa City, IA 52242-1419}
\date{} 
\date{ \today \\ preliminary version}
\maketitle

%%%%%%%%%%%%%%%%%%%%%%%%%%%%%%%%%%%%%%%%%%%%%%%%%%%

\section*{ \ \ \ \ \ \ Abstract}
\begin{quote}
{\small A classical theorem of Scheunert
on $G$-color Lie algebras, asserts  in the case of  finitely generated
abelian groups, one can twist the algebra structure and  the commutation 
bicharacter on $G$ by a 2-cocycle twist to a super-Lie $G$ graded, algebra.
In this paper we show that this can be done for an arbitrary group. 
}
\end{quote}

\section*{Introduction and notation}

We recall first the following definitions (see [Sch] and [Mo]):

Let $G$ be a group and 
$ \chi : G \times G \rightarrow k^{*} $
a bicharacter on $G$, i.e. a bimultiplicative morphism.
We assume that $ \chi $ is symmetric, i.e. $ \chi (h,g) \chi (g,h)=1 $ 
for all $h,g \in G$. 
Since $k$ is abelian it follows that $ \chi $ is trivial  for every 
commutator in $G$ so it factors through
 $ G^{ab} \times G^{ab}\rightarrow k^{*} $. 
Therefore from now on we assume $G$ to be abelian.

In this paper we assume that  $k$ is an algebraically 
closed field, $e$ will denote the the neutral element of $G$.

We call $L$ a $G$-color Lie algebra over $k$ with commutation factor
$\chi $ if $L$ is a $G$-graded 
$k$-vector space and the bracket $ [ \;, \; ] : L \times L \ra L$
satisfies:

$[\; a\;,\; b \; ] = - \chi (h,g) [\; b\;,\; a \; ]$

\no $  \chi (g,k) [\; a\;, [\; b\;,\; c \; ] \; ] +
\chi (k,h) [\; c\;, [\; a\;,\; b \; ] \; ] +
\chi (h,g) [\; b\;, [\; c\;,\; a \; ] \; ]$
for all $ a \in L_{g},  b \in L_{h},  c \in L_{k} $

We say that $\sigma$  is a $2$ cocycle on $G$ if 
$ \sigma  : G \times G \rightarrow k^{*} $ , satisfies

$ \sigma (a,bc)\sigma (b,c)=\sigma (a,b)\sigma (ab,c)$.
 Then we can define a new bracket
$ [ \;, \; ]_{\sigma}: L\times L \rightarrow L$
 by:

$[\; a\;,\; b \; ]_{\sigma } = \sigma (g,h) [\; a\;,\; b \; ]$
for all $ a \in L_{g},  b \in L_{h}$.

If $\sigma $ is a $2$ cocycle then
 $ \chi _{\sigma}(g,h)= \chi (g,h) \sigma (g,h) \sigma ^{-1}(h,g)$ 
is a bicharacter. We denote by $L^{\sigma }$ the (new) $G$-color 
 Lie algebra structure on $L$ for this new bracket
 $ [ \;, \; ]_{\sigma}$ and the commutation factor given by the 
twisted bicharacter  $ \chi _{\sigma} $. 

Let $  G_{+}=\{ g | \chi(g, g)=1 \}$, this is a subgroup of $G$
of index at most $2$, we call these the even elements in $G$. 
Define the odd elements by:  
$ G_{-}=\{ g | \chi(g, g)= -1\}$, then $G=  G_{+} \bigcup  G_{-}$.

We may define now $\chi _{o}( g | h )= 1$ iff at least
 one of $g$ or $h$ is even
else if both
$G$ and $H$ are odd $ \chi _{o}(g | h )= -1$

% State here Scheunert Thm.
Scheunert's theorem [Sch] shows that for a $G$ color Lie algebra $L$ with
bicharacter $\chi$ and  $G$  a finitely generated 
abelian group there exists a $2$-cocycle $\sigma $ on $G$ such that 
the bicharacter $ \chi _{\sigma} = \chi _{o}$. Thus $L^{\sigma}$
 can be regarded as a
  $\Z_{2}$ graded Lie algebra with  the $\Z_{2}$ (super)bicharacter 
$\chi _{o}$.

\section*{The proof for an arbitrary abelian group  $G$ }

%In fact besides Schuenert's paper there was one more result on twisting 
%To do that we need an ordered  system of generators of 
%$G$ s.t. any generator has at most one defining relation.

In this section we prove:
\bet
Let G be any abelian group and let $L$ be a $G$ color Lie algebra with\
 commutation factor $\chi $ Then there exists a bimultiplicative 
$2$-cocycle $\sigma $ on $G$ such that the twisted color Lie 
algebra $L^{\sigma}$ is  a super-Lie algebra with commutation factor 
$\chi _{o}$ .
\eet
Proof.

\non Like in the original proof in [Sche]  we may  change $\chi $ 
to $\chi \chi _{o}$
so that we may assume that $\chi (g,g)=1$ for all $g \in G$. 
We  show then that if  $\chi (g,g)=1$  for all $g \in G$ then 
there is a cocycle $\sigma $  on $G$ with 
$\chi (g,h) = \sigma (g,h) \sigma ^{-1} (h,g)$ for all $ g,h \in G$. 
Note that any bimultiplicative map is automatically a $2$-cocycle.

To do that we shall use Zorn's lemma.
Define a  family of subgroups of $G$ 

\no ${\cal F} = \{ ( H, \sigma _{H}) | H \mbox{subgroup of $G$,}
 \sigma _{H}) \mbox{ bilinear $2$-cocycle on $H$, } \chi (g,h) = 
\sigma (g,h) \sigma ^{-1} (h,g) \forall  g,h \in H \}$ 

We order this family by $( H', \sigma _{H'}) \preceq ( H'', \sigma _{H''})$
iff $H' \subseteq H''$ and ${\sigma _{H''}}|_{H'} = \sigma _{H'}$.
It is clear that $e \in \cal F$ so that $\cal F$ is non-void.
When in the sequel it  is clear on what subgroup $\sigma $ is  defined, we
 shall not show any more the indices. 

This way $\cal F $ is inductively ordered and by Zorn's lemma there exists a
maximal element of $\cal F $, say $( K, \sigma _{K})$.
Assume $K \neq G$. We shall prove this contradicts the maximality of $K$.

Let $t$ be an element in $G$ that does not belong to $K$.
We look at the subgroup $<t>$ generated by $t$.

If  $<t>\bigcap K = \{e\}$ then let $L=<t>\times K $. Define 
$\sigma (k, t) =  \chi ( k,t) $ and $\sigma ( t,k) = 1$ for all $k \in K$
and extend $ \sigma $ bimultiplicatively. Since there are no new relations 
this is well defined  and one can see that 
 $\chi (g,h) = \sigma (g,h) \sigma^{-1} (h,g)$ holds on $L$.

If $<t>\bigcap K = <t^{n}>$ then there are 
some $ k_{1}, k_{2} \ldots k_{m} \in K$  and some positive integers
 $ n_{1}, n_{2} \ldots n_{m}$ such that
$ t^{n} =  k_{1} ^{n_{1}} k_{2} ^{n_{2}} \ldots k_{m}^{n_{m}}$.
More than one such relations is possible but we just select one, say with 
a minimal $m$.

Define now $L=<t,K>$ to be the subgroup generated by $K$ and $t$.
We need to extend  $\sigma $ to $L$ in such a way that :

 1) $\sigma $ is well defined and  bimultiplicative on $L$

 2) $\chi (g,h) = \sigma (g,h) \sigma ^{-1} (h,g)$ for all $ g,h \in L$. \\

\no Because  $ t^{n} =  k_{1} ^{n_{1}} k_{2} ^{n_{2}} \ldots k_{m}^{n_{m}}$
it is clear that for any $u \in L$ we have
$ \sigma (u,t^{n}) = \sigma (u, k_{1} ^{n_{1}} k_{2} ^{n_{2}} 
\ldots k_{m}^{n_{m}}  )$.

This means is $ \sigma (u,t^{n})$ is already determined, so loosely 
speaking we may say 
$\sigma (u, t) = {}^n \! \! \sqrt{{ \prod_{i} \sigma ( u, k_{i}^{n_{i}})}} $

The problem is that while we  have $n$-th roots, $k$ being algebraically
 closed,  we do not have a uniform {\em radical function} (say like the real 
radical), so we need to make sure that we define $\sigma $ as a  function
multiplicative on both first and second variable.

\no We start by  defining a multiplicative function
 $f(u)=\sigma (u,t), f: K \rightarrow k^{*}$ (multiplicative in $u$) 
such that:

 $ f(u)^n =  \sigma (u,t^{n}) = 
\sigma (u, k_{1} ^{n_{1}} k_{2} ^{n_{2}} \ldots k_{m}^{n_{m}} )$ 
and $f(t^n)=1$

We let $\cal M $ be the family of subgroups of $K$ that contain $<t^{n}>$, 
on which $f$ can  be defined with the above properties, ordered by set 
inclusion and by the requirement that  $f$ extends  from the small 
subgroup to the bigger one. 

Then $\cal M $ is non-void and inductively ordered hence 
it has a maximal element $M$. If this maximal element is not $K$  itself
say $M \subset K $ and $ M \neq K $ then we may contradict the maximality 
of $M$. 

For an $w\in K \; - \; M$ we extend  $f$ to $<w, M>$ by: 

If $ <w> \bigcap M =\{e\} $ then let $f(w)$ be any selection of 
${}^n \sqrt{ \prod_{i} \sigma ( w, k_{i}^{n_{i}})} $.
This works  since 
$<w,M>=<w> \times M$ and  contradicts the maximality of $M$ unless $M=K$.

Else if $ <w> \bigcap M =<w^{r}> $ and $w^{r}= \prod z_{j}^{r_{j}}$, 
(a finite product)  define:

$f(w)={}^{rn}\! \! \sqrt{\prod_{i,j} {\sigma ( z_{j}, k_{i})}^{n_{i} r_{j}}}$

This contradicts again the maximality of $M$ and it means there is 
a multiplicative mapping  

$f(u)=\sigma (u,t): K \rightarrow k^{*}$ such that
 $ f(u)^{n} =  \sigma (u,t^{n}) = \sigma (u, k_{1} ^{n_{1}} k_{2} ^{n_{2}}
 \ldots k_{m} ^{n_{m}})$ and $f(t^{n})=1$

We use now the required relation to move $u$ on the right side by defining 
now
an analog of $f$ on the ``right'': 

$\sigma (t, u) = \chi(t,u) \sigma (u,t)$ 

This is  multiplicative in the second variable, i.e. in $u$, because $\chi $
is bimultiplicative and also $f$ is multiplicative.

Since $\chi(g,g)=1$ was granted  we define $\sigma (t,t)=1$, this is 
consistent with the previous definitions (and this was the reason 
we asked $f(<t^n>)=1$).

Now we define $\sigma $ on all $<t,K>$ by

$$ \sigma (t^{\alpha } u,t^{\beta } v) = {\sigma (t,v)}^{\alpha} 
 {\sigma (u,t)}^{\beta }  \sigma (u,v)$$

This is bimultiplicative because of the way it was defined.
One can use the fact that $f$ is multiplicative, to show that $\sigma $ 
is   well defined. One needs to show that $\sigma $
 respects the relation: 
 $ t^{n} =  k_{1} ^{n_{1}} k_{2} ^{n_{2}} \ldots k_{m}^{n_{m}}$, 
when substituted on either side. For this we look at reduced forms 
of $t^{\alpha } u$ , 
with $\alpha < n$. The relation holds because of the way that $f(t)$ was
defined.
This way we contradict now the maximality of $K$ so we may conclude $K=G$
and the proof of  our theorem.
%\qed

\brk
There is another really interesting instance of  twisting in the
paper by Artin-Schelter-Tate [AST]. It is proved there that
the multiparametric quantum general linear group is a twist 
of the standard quantization of the general linear group.

In fact we are interested in the result of Proposition 1 in [AST], 
where $G$ is a free abelian 
group of dimension $n<\infty$, it is  shown that any cocycle 
cohomology class in $ H^2(G,k^{*})$ contains exactly one 
bicharacter on $G$. We conjecture this is the case for an arbitrary 
abelian group somehow along a construction similar to that of $\sigma $
in the proof above.
\erk

\brk
In fact the proof here does not fully use the fact that $k$ is
 algebraically closed. Assume that we use transfinite induction to find 
the following presentation for $G$: $G$ is given by
a system of generators $\{ t_{\lambda} \}_{\lambda\in \Lambda} $ 
such that for each generator $t_{\lambda}$ there is a 
unique relation $r(t_{\lambda}):  
t^{n_{\lambda}} =  k_{1} ^{n_{1}} k_{2} ^{n_{2}} \ldots k_{m}^{n_{m}}$.

In our proof we  only used the fact that  $k$ was closed
under radicals of orders equal to the numbers $n_{\lambda}$ 

\erk

\no {\bf Corolaries}. The ones given in [Sch]: PBW bases and Ado's Theorem.

\end{document}